\documentclass{amsart}
\usepackage{amssymb, amsmath, amsfonts, amscd}
\usepackage{hyperref}

\input xypic

\theoremstyle{plain}
\newtheorem{theorem}{Theorem}[section]
\newtheorem{corollary}[theorem]{Corollary}

\newtheorem{proposition}[theorem]{Proposition}

\newtheorem{example}[theorem]{Example}

\theoremstyle{definition}
\newtheorem{definition}[theorem]{Definition}

\theoremstyle{remark}

\numberwithin{equation}{theorem}

\renewcommand{\O}{\mathcal{O} }

\newcommand{\Spec}{\operatorname{Spec} }
\renewcommand{\P}{\operatorname{P} }
\renewcommand{\H}{\operatorname{H} }

\renewcommand{\P}{\mathbb{P} }
\newcommand{\Proj}{\operatorname{Proj}}

\begin{document}

\title{On the birational invariance of the arithmetic genus and Euler characteristic}



\author{Helge \"{O}ystein Maakestad \\ email \href{mailto:h_maakestad@hotmail.com}{\text{h\_maakestad@hotmail.com}} }

\email{\text{h\_maakestad@hotmail.com}}

\address{Address: Tempelveien 112, 3475 Saetre i Hurum, Norway}

\keywords{projective variety, hypersurface, arithmetic genus, geometric genus, Euler characteristic, birational invariance}

\thanks{}

\subjclass{}

\date{May 2019}

\begin{abstract} The aim of this note is to use elementary methods to give a large class of examples of projective varieties $ Y \subseteq \mathbb{P}^d_k$ over a field $k$ with the property that
$Y$ is not isomorphic to a hypersurface $H\subseteq \mathbb{P}^N_k$ in projective space $\mathbb{P}^N_k$ with $N:=dim(Y)+1$. We apply this construction to the study of the arithmetic genus
$p_a(Y)$ of $Y$ and the problem of determining if $p_a(Y)$ is a birational invariant of $Y$ in general.  We give positive dimensional families of pairs of projective varieties $(Y, Y')$ in any dimension 
 $d \geq 4$  where $Y$ is birational to $Y'$, but where $p_a(Y)\neq p_a(Y')$. The arithmetic genus is known to be a birational invariant in any dimension $d$ over an algebraically closed field of characteristic zero. The varieties $Y,Y'$ may be singular. We get similar results for the Euler characteristic $\chi(Y)$.
\end{abstract}

\maketitle

\tableofcontents

\section{Introduction} The problem of defining birational invariants for smooth projective varieties has a long and complicated history in algebraic geometry. The geometric genus $p_g(Y)$ and arithmetic genus $p_a(Y)$
of a projective variety  $Y$ has been around since the beginning of the study of projective geometry and projective varieties.
It is well known that the geometric genus $p_g(Y)$ of a smooth projective variety $Y$ over an algebraically closed field $k$ is a birational invariant. The arithmetic genus $p_a(Y)$  of a smooth projective variety $Y$ over 
an algebraically closed field of characteristic zero is by Hodge theory a birational invariant. In this paper we give positive dimensional families of pairs of birational projective varieties $(Y,Y')$ in any dimension
$d\geq 4$ with $p_a(Y)\neq p_a(Y')$. The varieties $Y,Y'$ may be singular and the base field $k$ is arbitrary. Hence the arithmetic genus $p_a(-)$ is not a birational invariant in general. We get a similar result for the Euler characteristic.


\section{Projective varieties that are not hypersurfaces}

The aim of this section is to use elementary methods to give a large class of examples of projective varieties $Y$ of finite type over a field $k$ with negative arithmetic genus $p_a(Y)$, and to use this class of examples 
to study the birational invariance of the arithmetic genus $p_a(Y)$ and Euler characteristic $\chi(Y)$.

In Theorem \ref{prod} we construct a projective variety 
$Y_{d,l}:=\H_d \times H_l$ as a product of two hypersurfaces $H_d$ and $H_l$ of degrees $d$ and $l$ with the property that $p_a(H_d)=0$ and $dim(H_l)$ is odd. It follows $p_a(Y_{d,l})<0$. If a projective variety $Y$ has negative arithmetic genus it cannot be isomorphic to a hypersurface $H$ in a projective space $\P^n_k$ with $n=dim(Y)+1$, since a hypersurface $H$ always have $p_a(H)\geq 0$. Using this construction we give positive dimensional families of pairs $(Y,Y')$ of projective varieties in any dimension $d\geq 4$ where $Y$ and $Y'$ are birational but where $p_a(Y)\neq p_a(Y')$. The arithmetic genus is known to be a birational invariant for smooth projective varieties over an algebraically closed field of characteristic zero in any dimension (see \cite{griffiths}, page 494). We get a similar result for the Euler characteristic.

Let in this section $k$ be an arbitrary field and let all varieties be of finite type over $k$ as defined in \cite{hartshorne}, Chapter I.
Let $Y\subseteq \P^n_k$ be a projetive variety of dimension $r$ with Hilbert polynomial $P_Y(t)$. 

Note: The variety $Y$ is defined using a homogeneous prime ideal 
$I(Y) \subseteq k[x_0,..,x_n]$ and $Y\cong \Proj(S)$ where $S:=k[x_0,..,x_n]/I(Y)$. The Hilbert polynomial $P_Y(t)$ is defined using the graded ring $S$ as done in \cite{hartshorne} I.7.5, 
hence the polynomial $P_Y(t)$ depends on the ring $S$ and the embedding $Y\subseteq \P^n_k$. We use the polynomial $P_Y(t)$ to define the arithmetic genus $p_a(Y)$.

\begin{definition} Let $p_a(Y):=(-1)^r(P_Y(0)-1)$ be the  \emph{arithmetic genus} of $Y$.
\end{definition}

The arithmetic genus $p_a(Y)$ is independent of choice of embedding $Y \subseteq \P^n_k$ and is an important invariant of the variety $Y$ (see \cite{hirzebruch}). 
We may in the case when $k$ is algebraically closed use the Euler characteristic $\chi(\O_Y)$ of the structure sheaf $\O_Y$ to define it: 

\begin{proposition} \label{harteuler} Let $Y\subseteq \P^n_k$ be a closed subvariety where $k$ is an algebraically closed field. The following formula holds:
\[ p_a(Y)=(-1)^r(\chi(\O_Y)-1) \]
\end{proposition}
\begin{proof}
The proof follows from \cite{hartshorne}, Exercise III.5.3.
\end{proof}

The structure sheaf $\O_Y$ is independent of choice of closed embedding $Y\subseteq \P^n_k$. From this it follows $p_a(Y)$ is independent of choice of embedding into a projective space.
The calculations in the paper are based on the  following formulas for $p_a$: 

\begin{proposition} \label{hart} Let $H_d \subseteq \P^N_k$ be a hypersurface of degree $d$ and let $Y\subseteq \P^m_k$ and $Z\subseteq \P^n_k$ be algebraic varieties of dimension
$r$ and $s$. The following formulas hold for the arithmetic genus $p_a(-)$.
\begin{align}
&\label{hyper}p_a(H_d)=\binom{d-1}{N} \\
&\label{product}p_a(Y\times Z)=p_a(Y)p_a(Z)+(-1)^s p_a(Y)+(-1)^rp_a(Z).
\end{align}
\end{proposition}
\begin{proof} This is Exercise I.7.2 in \cite{hartshorne}.
\end{proof}

Note: When $d-1<N$ we define $\binom{d-1}{N}:=0$. Hence if $H_d \subseteq \P^N_k$ is a hypersurface of degree $d$ with $d-1<N$ it follows $p_a(H_d)=0$. In particular since we may realize $\P^n_k$ as a hypersurface $H_1\subseteq \P^{n+1}_k$ of degree one, it follows $p_a(\P^n_k)=0$ for all $n\geq 1$.

\begin{theorem}\label{prod}  Assume $H_d \subseteq \P^{2n}_k$ be a hypersurface of degree $d$ with $d-1<n$ and let $H_l\subseteq \P^m_k$ be a hypersurface of degree $l$.
it follows
\[ p_a(H_d \times H_l)=-\binom{l-1}{m}<0.\]
\end{theorem}
\begin{proof} From Proposition \ref{hart} we get since $p_a(H_d)=0$ the formula
\[ p_a(H_d\times H_l)=p_a(H_d)p_a(H_l)+(-1)^{m-1}p_a(H_d)+(-1)^{2n-1}p_a(H_l)=-\binom{l-1}{m} .\]
Hence $p_a(H_d \times H_l)<0$ in general: Choose an integer $l$ with $l-1\geq m$.
\end{proof}

\begin{corollary} Let $H_d,H_l$ be hypersurfaces satisfying the conditions in Theorem \ref{prod}. It follows the product variety $H_d \times H_l$ cannot be isomorphic to a hypersurface $H_e\subseteq \P^N_k$
of degree $e$ in projective $N$-space with $N:=dim(H_d \times H_l)+1$.
\end{corollary}
\begin{proof} Since $p_a(H_e)=\binom{e-1}{N} \geq 0 $ and $p_a(H_d\times H_l)<0$ there can be no isomorphism $H_d \times H_l \cong H_e$ since $p_a(-)$ is invariant under isomorphism.
\end{proof}

\begin{definition}
Two varieties $Y,Y'$ of the same dimension $d$ are \emph{birational} if there are open subvarieties $U\subseteq Y$ and $U' \subseteq Y'$ and an isomorphism of varieties $U\cong U'$.
\end{definition}

\begin{example} Affine $n$-space and projective $n$-space are birational.\end{example}

Let $\mathbb{A}^n_k:=\Spec(k[y_1,\ldots , y_n])$ and $\P^n_k:=\operatorname{Proj}(k[x_0,\ldots ,x_n])$. It follows 
there is an open embedding $\phi: \mathbb{A}^n_k \rightarrow \P^n_k$, hence $\mathbb{A}^n_k$ is birational to $\P^n_k$.

\begin{corollary} The arithmetic genus is not a birational invariant in general.
\end{corollary}
\begin{proof} Let $H_d,H_l$ be hypersurfaces satisfying the conditions in Theorem \ref{prod}. It follows $p_a(H_d \times H_l) <0$. Any projective variety $X\subseteq \P^N_k$ of dimension $r$ is by \cite{hartshorne}, Proposition I.4.9 birational
to a hypersurface $H_e \subseteq \P^{r+1}_k$ of degree $s$ for some integer $s$. If the arithmetic genus was a birational invariant in general, we would have
\[ p_a(H_d \times H_l) = p_a(H_e) \geq 0 .\]
 This gives a contradiction since $p_a(H_d \times H_l)<0$, and the Corollary follows.
\end{proof}

\begin{corollary} \label{maincorr}Let $\P^1_k$ be the projective line and let $H_l\subseteq \P^n_k$ be a hypersurface of degree l, with $l-1\geq n$ and $n\geq 4$. It follows $dim(\P^1_k\times H_l)=n$.
There is an equality $p_a(\P^1_n \times H_l)=-\binom{l-1}{n}<0$. There is a birational isomorphism $\P^1_k \times H_l \cong H_e$ where $H_e \subseteq \P^{n+1}_k$ is a hypersurface of degree $e$
and
\[ p_a(\P^1_k \times H_l) \neq p_a(H_e)=\binom{e-1}{n+1}.\]
Hence we get examples of pairs of birational varieties $(\P^1_l\times H_l, H_e)$ with different arithmetic genus in any dimension $dim(\P^1_k \times H_l) \geq 4$.
\end{corollary}
\begin{proof} We get
\[ p_a(\P^1_k \times H_l)=p_a(\P^1_k)p_a(H_l)+(-1)^{n-1}p_a(\P^1_k)+(-1)^1p_a(H_l)=-p_a(H_l)=-\binom{l-1}{n}\]
since $p_a(\P^1_k)=0$. There is by \cite{hartshorne}, Proposition I.4.9 a hypersurface $H_e\subseteq \P^{n+1}_k$ and a birational morphism $\P^1_k\times H_l \cong H_e$. Since
$p_a(\P^1_k \times H_l) <0$ and $p_a(H_e) \geq 0$ it follows $p_a(\P^1_k \times H_l) \neq p_a(H_e)$.  Since $dim(\P^1_k \times H_l)=n\geq 4 $ the Corollary follows. 
\end{proof}

\begin{corollary} \label{eulerchar} Let $k$ be an arbitrary algebraically closed field and let $H_l\subseteq \P^n_k$ be a hypersurface of degree $l$ with $l-1\geq n$ and $n \geq 4$. There is a birational isomorphism 
$\P^1_k \times H_l \cong H_e$ where $H_e \subseteq \P^{n+1}_k$ is a hypersurface of degree $e$ and $\chi(\P^1_k \times H_l) \neq \chi(H_e)$. Hence the Euler characteristic is not a birational invariant in general.
\end{corollary}
\begin{proof} By Proposition \ref{harteuler} it follows $p_a(X)=(-1)^r(\chi(X)-1)$ with $r:=dim(X)$. Hence since $p_a(\P^1_k \times H_l) \neq p_a(H_e)$ it follows $\chi(\P^1_k \times H_l) \neq \chi(H_e)$.
\end{proof}

Let the hypotheses be as in Corollary \ref{eulerchar}.
By the proof of Corollary \ref{maincorr} we get positive dimensional families of pairs of birational projective varieties $(Y,Y')$ over $k$ with $dim(Y)=dim(Y') \geq 4 $ and where $\chi(Y) \neq \chi(Y')$.

\begin{example} \label{arbitrary}Varieties with arbitrary large arithmetic genus. \end{example}

For a projective variety $X\subseteq \P^n_k$ we get the following result: Let $k$ be an arbitrary field and let $X$ be a projective varitety with $p_a(X)\neq 0$. Assume $p_a(X)<0$ and let $H_e:=Z(f) \subseteq \P^{n+1}_k$ be a hypersurface of degree $e$ with $X\cong H_e$ birational. Let $g\in \H^0(\P^{n+1}_k,\O(d))$ be a homogeneous polynonmial of degree $d$ and let $H_{e+d}:=Z(fg)$. 

\begin{theorem} It follows there is a birational isomorphism
$X\cong H_{e+d}$ and $p_a(H_{e+d})=\binom{e+d-1}{n+1}$. The polynomial $g$ is arbitrary and the arithmetic genus $p_a(H_{e+d})$ can be arbitrary large. If $p_a(X)$ is positive
we get a similar result for $Y:=\P^1_k \times X$.
\end{theorem}
\begin{proof} The proof is similar to the proof of Corollary \ref{maincorr}.
\end{proof}

Hence we get for any projective variety $X$ over $k$ with non-zero arithmetic genus, a family $(Y,H_{e+d})$ of pairs of birational varieties, but where the arithmetic genus of $H_{e+d}$ can be made arbitrary large.
A similar result holds for the Euler characteristic when the base field is algebraically closed. 
Hence it is easy to give examples of pairs of projective varieties $(Y,Y')$ that are birational, but where the arithmetic genus and Euler characteristic differ.

\begin{example} Birational invariance of the arithmetic genus, geometric genus and Euler characteristic. \end{example}

The arithmetic genus $p_a(Y)$ of a projective variety has been around since the beginning of the study of projective geometry and projective varieties.
You may find some information on the problem of deciding if the arithmetic genus is a birational invariant in \cite{hartshorne}. 
In \cite{hartshorne} Exercise III.5.3 it is proved that the arithmetic genus of a non-singular projective curve $C$
over an algebraically closed field  is a birational invariant.  In \cite{hartshorne}, Corollary V.5.6 it is proved that the arithmetic genus
of a nonsingular projective surface $S$ over an algebraically closed field is a birational invariant.  In \cite{griffiths}, page 494 you find a proof that the Euler characteristic $\chi(X)$ 
of a smooth projective variety $X$ of dimension $r$ over the complex numbers is a birational invariant. It follows the arithmetic genus $p_a(X)=(-1)^r(\chi(X)-1)$ is a birational invariant.
In Corollary \ref{maincorr} and Corollary \ref{eulerchar} we get  positive dimensional families of pairs of birational projective varieties $(Y,Y')$ over any field $k$ 
with $p_a(Y)\neq p_a(')$ and $\chi(Y)\neq \chi(Y')$  in any dimension $d\geq 4$. The varieties $Y,Y'$ may be singular.

\textbf{Acknowledgements.} Thanks to J. B. Bost for some information on Hodge theory and the arithmetic genus.

\end{document}